\documentclass[10pt]{article}
\usepackage{enumerate,amsmath,xspace,amssymb,mathrsfs,latexsym,xy}
\xyoption{all}\xyoption{2cell}
\UseAllTwocells
\CompileMatrices
\input{diagxy}

\let\named\label  

\newenvironment{proof_env}{\begin{trivlist}\item[]{Proof\/}.%
   \hspace{1em}}{\parfillskip=0pt\hfill%
   $\scriptstyle \Box$\end{trivlist}}

\newcommand{\proof}{\begin{proof_env}}
\newcommand{\qed}{\end{proof_env}}

\newcommand{\B}{{\ensuremath{\mathscr{B}}}\xspace}

\newcommand{\bfS}{\ensuremath{\textbf{S}}}

\newcommand{\op}{\ensuremath{^{\textnormal{op}}}}

\newcommand{\<}{\langle}
\renewcommand{\>}{\rangle}

\renewcommand{\phi}{\varphi}

\def\mto{\mapsto}

\def\ra{\xy\morphism(0,0)|a|/>/<150,0>[`;]\endxy}
\def\la{\xy\morphism(0,0)|a|/<-/<150,0>[`;]\endxy}
\def\pro{\xy\morphism(0,0)|a|/>/<150,0>[`;]\endxy}
\def\arr#1{\xy\morphism(0,0)|a|/>/<150,0>[`;#1]\endxy}

\def\epi{\xy\morphism(0,0)|a|/>>/<150,0>[`;]\endxy}

\def\lra{\xy\morphism(0,0)|a|/>/<200,0>[`;]\endxy}
\def\arrow#1{\xy\morphism(0,0)|a|/>/<200,0>[`;#1]\endxy}
\def\2cell{\xy\morphism(0,0)|a|/=>/<100,0>[`;]\endxy}
\def\rev2cell{\xy\morphism(0,0)|a|/<=/<100,0>[`;]\endxy}

\def\adj{\dashv}

\def\iso{\cong}
\def\equ{\:\simeq\:}
\def\id{\rm id}

\def\colimit#1{\xy\morphism(0,-40)|a|/>/<120,0>[`;{\rm lim}]%
 \place(50,-80)[{\scriptstyle #1}]\endxy}

\newbox\anglebox
\setbox\anglebox=\hbox{\xy \POS(75,0)\ar@{-} (0,0) \ar@{-} (75,75)\endxy}%

\def\Se{\mbox{\sl Set}}

\def\PSh{\mbox{\sl PSh}}

\def\bA{\mathbb{A}}
\def\bB{\mathbb{B}}
\def\bC{\mathbb{C}}
\def\bD{\mathbb{D}}

\def\bP{\mathbb{P}}

\def\bS{\mathbb{S}}
\def\bU{\mathbb{U}}

\def\bX{\mathbb{X}}

\newtheorem{theorem}{Theorem}[section]
\newtheorem{proposition}[theorem]{Proposition}
\newtheorem{corollary}[theorem]{Corollary}
\newtheorem{lemma}[theorem]{Lemma}

\newtheorem{remark}[theorem]{Remark}
\newtheorem{Claim}{Claim}
\newcommand{\dom}{\mathbf{d}}
\newcommand{\ran}{\mathbf{r}}

\title{Characterizations of Morita equivalent inverse semigroups}
\author{J.~Funk
\and M.~V.~Lawson
\and B.~Steinberg}

\begin{document}
\maketitle

\begin{abstract}

We prove that four different notions of Morita equivalence for inverse semigroups motivated by, respectively, $C^{\ast}$-algebra theory,
topos theory, semigroup theory and the theory of ordered groupoids are equivalent.
We also show that the category of unitary actions of an inverse semigroup
is monadic over the category of \'etale actions.
Consequently, the category of unitary actions of an inverse semigroup
is equivalent to the category of presheaves on its Cauchy completion.
More generally, we prove that the same is true for the category of closed actions,
which is used to define the Morita theory in semigroup theory,
of any semigroup with right local units.\\

\end{abstract}

\section{Introduction}\named{intro}\setcounter{theorem}{0}

The Morita theory of unital rings was introduced by Morita
in 1958 \cite{Morita}:
two such rings are {\em Morita equivalent\/}
if their categories of left modules are equivalent.
This definition provides a classification of rings
that is weaker than isomorphism but still useful; in particular,
the Artin-Wedderburn theorem can be interpreted in terms of Morita equivalence.
There are at least two important characterizations of Morita equivalence.
The first uses the notion of invertible bimodules \cite{B}:
rings $R$ and $S$ are Morita equivalent if and only if there is an $(R,S)$-bimodule $X$ and an $(S,R)$-bimodule $Y$ such that
$X \otimes Y \cong R$ and $Y \otimes X \cong S$.
The second uses rings of matrices and full idempotents \cite{Lam}:
rings $R$ and $S$ are Morita equivalent if and only if $R$
is isomorphic to a ring of the form $eM_{n}(S)e$
where $e$ is a full idempotent meaning that $M_{n}(S) = M_{n}(S)eM_{n}(S)$.
These results have been the model for analogous definitions made for other structures: for example, monoids \cite{B,K} and (small) categories \cite{EZ}.
The theory has also been extended to classes of non-unital rings \cite{A,AM}.
This in turn inspired a Morita theory for semigroups \cite{T1,T2,T3} due to Talwar.

This paper concerns the Morita theory of a class of semigroups
called inverse semigroups.
These are one of the most interesting classes of semigroups
with connections to diverse branches of mathematics.
They are the abstract counterparts of pseudogroups of transformations and can be viewed as carriers of information about partial symmetries \cite{L2}.
There are also very close connections between inverse semigroups and topoi \cite{F,F1,L4}.
We define them as follows.
A semigroup $S$ is {\em (von Neumann) regular\/} if for each $s \in S$ there exists $t \in S$, called an {\em inverse} of $s$,
such that $s = sts$ and $t = tst$.
If each element of a regular semigroup has a unique inverse,
then the semigroup is said to be {\em inverse}.
We denote the unique inverse of an element $s$ in an inverse semigroup by $s^{\ast}$ in this paper.
Equivalently, a regular semigroup $S$ is inverse if its sets
of idempotents $E(S)$ forms a commutative subsemigroup.
The set of idempotents $E(S)$ of an inverse semigroup is ordered
when we define $e \leq f$ whenever $e = ef = fe$.
With respect to this order, the set $E(S)$ is a meet-semilattice
in which $e \wedge f = ef$.
For this reason, the set of idempotents of an inverse semigroup is usually referred to as its {\em semilattice of idempotents}.

Let us make some definitions for arbitrary semigroups.
Let $X$ be a set and $S$ a semigroup.
We say that $X$ is a {\em right $S$-set} if there is a function
$X \times S \ra X$, given by $(x,s) \mapsto xs$,
such that $x(st) = (xs)t$ for all $x \in X$ and $s,t \in S$.
Left $S$-sets are defined dually.
If $S$ and $T$ are both semigroups that act on the set $X$
on the left and right respectively
in such a way that $(sx)t = s(xt)$ for all $s \in S$, $t \in T$ and $x \in X$
then we say that $X$ is an {\em $(S,T)$-biset}.
In this paper, we shall usually only deal with right $S$-sets,
so that we shall usually omit the word `right' in what follows.
An $S$-set $X$ is said to be {\em unitary\/} if for every $x\in X$
there are $s\in S$, $y\in X$ such that $ys=x$.
We write $XS=X$.

This paper is motivated by the fact that there are no fewer than four
possible definitions of Morita equivalence for inverse semigroups:
\begin{enumerate}

\item strong Morita equivalence;

\item topos Morita equivalence;

\item semigroup Morita equivalence;

\item enlargement Morita equivalence.

\end{enumerate}
We shall now define each of these notions.\\

\noindent
{\bf 1. Strong Morita equivalence}\\

Inverse semigroups $S$ and $T$ are said to be {\em strongly Morita equivalent} \cite{S} if there is an {\em equivalence biset\/} for $S$ and $T$;
by definition, this consists of a set $X$,
which is an $(S,T)$-biset equipped with surjective functions
$$
\langle -,- \rangle \colon \: X \times X \ra S\;,
\text{ and }
[-,-] \colon \:X \times X \ra T
$$
such that the following axioms hold,
where $x,y,z \in X$, $s \in S$, and $t \in T$:
\begin{description}

\item[{\rm (M1)}] $\langle sx,y \rangle = s\langle x,y \rangle$

\item[{\rm (M2)}] $\langle y,x \rangle = \langle x,y \rangle^{\ast}$

\item[{\rm (M3)}] $\langle x,x \rangle x = x$

\item[{\rm (M4)}] $[x,yt] = [x,y]t$

\item[{\rm (M5)}] $[x,y] = [y,x]^{\ast}$

\item[{\rm (M6)}] $x[x,x] = x$

\item[{\rm (M7)}] $\langle x,y \rangle z = x [y,z]\;$.

\end{description}
This definition is motivated by Rieffel's notion
of an equivalence bimodule \cite{S},
and is well adapted to the natural affiliation of inverse semigroups with both \'etale topological groupoids and $C^{\ast}$-algebras \cite{P};
in particular,

\begin{itemize}
\item
if $S$ and $T$ are strongly Morita equivalent, then their associated \'etale groupoids \cite{P} are Morita equivalent;

\item
if $S$ and $T$ are strongly Morita equivalent,
then their universal and reduced $C^{\ast}$-algebras are strongly Morita equivalent \cite{S}.
\end{itemize}

\noindent
{\bf 2.~Topos Morita equivalence }\\

Whereas strong Morita equivalence takes the bimodule aspect of classical Morita theory as it starting point, another natural starting point is actions.
Let  $S$ be an inverse semigroup.
Then $S$ acts on its semilattice of idempotents $E(S)$ when we define
$e \cdot s = s^*es$.
We call this the {\em Munn $S$-set}.
An $S$-set $X$ paired with an $S$-set map
$X \arr{p} E(S)$ to the Munn $S$-set,
such that $x\cdot p(x) = x$,
is called an {\em \'etale right $S$-set\/} \cite{F1}.
We denote the category of \'etale right $S$-sets by {\sl \'Etale\/}.
The category {\sl \'Etale\/} is a topos, sometimes called
the classifying topos of  $S$ and is also denoted by $\B(S)$.\footnote{
The term `classifying topos' and its $\B$ notation more generally refer to
the topos associated with an \'etale, or even localic, groupoid \cite{MM}.
An ordered groupoid is \'etale in this sense.
It is not difficult to see that the definition $\B(S)=\B(G(S))$
ultimately amounts to the category of \'etale $S$-sets.}
{\sl \'Etale\/} is in a sense the `space' of $S$, but the following
`categorical' description of it is sometimes important for calculations.
With the inverse semigroup $S$, we may associate a left cancellative category
$$
L(S) = \{(e,s) \in E(S) \times S \colon \: es=s\}\;,
$$
whose composition is given by
$(e,s)(f,t) = (e,st)$, provided $s^{\ast}s = f$.
The objects of $L(S)$ can be identified with $E(S)$ and the arrow $(e,s)$ goes from $s^*s$ to $e$.  The identity at $e$ is $(e,e)$.
The category {\sl \'Etale\/} is equivalent to the category $\PSh(L(S))$ of presheaves on $L(S)$,
where a {\em presheaf} on a category is a contravariant
functor to the category of sets.
This result, which is used in \cite{F,F1,L4},
is essentially due to Loganathan \cite{L}.

We say that two inverse semigroups $S$ and $T$ are {\em topos Morita equivalent} if the categories  $\B(S)$ and $\B(T)$ are equivalent.
Steinberg \cite{S} proves that strong Morita equivalence
implies topos Morita equivalence, but whether the converse is true was left open.
We shall see later that they are indeed equivalent.\\

\noindent
{\bf 3.~Semigroup Morita equivalence}\\

The previous definition viewed inverse semigroups
within the context of topos theory.
They can of course be viewed simply as semigroups, and for a wide class of semigroups there is another definition of Morita equivalence.
Let $S$ be a semigroup with set of idempotent $E(S)$.
We say that $S$ has {\em right local units\/} if $SE(S)=S$.
Having {\em left local units} is defined dually and one says that
$S$ has {\em local units} if it has both left and right local units.
Inverse semigroups and more generally regular semigroups have local units.
We shall assume that $S$ is a semigroup with right local units.
Let $X$ be a set equipped with a right action $\mu\colon X\times S\ra X$.
The universal property of the tensor product yields an induced map $\mu\colon X\otimes_S S\ra X$ given by $x\otimes s\mto xs$.
Notice that $\mu$ is surjective precisely when the action is unitary.
One says that $X$ is {\em closed\/} if $\mu$ is also injective.
The category of closed $S$-sets will be denoted $S\text{-}\Se\/$.
Following Lawson and Talwar~\cite{AL,T1,T2,T3},
we say that two semigroups $S$ and $T$ with right local units
are {\em semigroup Morita equivalent\/} if the categories $S\text{-}\Se\/$ and $T\text{-}\Se\/$ are equivalent.
Talwar~\cite{T1} proves that if $S$ is an inverse semigroup,
then the closed right $S$-sets are precisely the unitary ones.
Thus, when $S$ is inverse $S\text{-}\Se\/$ is
the category of unitary right $S$-sets.

In the theory of semigroup Morita equivalence
another category plays an important role.
Let $S$ be any semigroup.
Then
$$
C(S) = \{(e,s,f) \in E(S) \times S \times E(S) \colon \: esf = s \},
$$
with the obvious partial binary operation, is a category called the {\em Cauchy completion\/} of $S$ (other terminology includes the idempotent
splitting and the Karoubi envelope).
The objects of $C(S)$ are again the idempotents of $S$.
A morphism $(e,s,f)$ of $C(S)$ may also be depicted $f\arr{s}e$.
In the case where $S$ is inverse, the category $L(S)$ is a subcategory of $C(S)$, although not necessarily full.
One identifies the arrow $(e,s)$ of $L(S)$ with $(e,s,s^*s)$.\\

\noindent
{\bf 4.~Enlargement Morita equivalence}\\

An inverse semigroup $S$ can also be regarded as a special kind of ordered groupoid $G(S)$ called an {\em inductive groupoid}.
An {\em ordered groupoid\/} $G$ is a groupoid internal to
the category of posets such that the domain map is a discrete fibration.
Equivalently, $G$ is an ordered groupoid if it is \'etale,
when regarded as a continuous groupoid
with respect to its downset (Alexandrov) topology~\cite{F,L2}.
The underlying set of $G(S)$ is $S$,
the groupoid product is the restricted product,
and the order is the natural partial order on $S$.
In this way, the category of inverse semigroups
can be embedded in the category of ordered groupoids.
We denote by $\dom$ and $\ran$ the domain and
range of an element of an ordered groupoid.
If $g$ and $h$ are elements of an ordered groupoid
such that $e = \dom (g) \wedge \ran (h)$ exists,
then we may define their {\em pseudoproduct}
by $g \circ h = (g \mid e)(e \mid h)$.
We refer the reader to \cite{L2} for the definitions and the basic theory.

We may extend some of the definitions
we have made earlier to classes of ordered groupoids.
Let $G$ be an arbitrary ordered groupoid.
We define the category $L(G)$ to consist of ordered pairs $(e,g)$,
where $\ran (g) \leq e$,
with product given by $(e,g)(f,h) = (e, g \circ h)$ when $\dom (g) = f$.
Observe that the pseudoproduct is defined.
This directly extends the definition
we made of this category in the inverse semgiroup case.
The {\em classifying topos\/} $\B(G)$ is by definition
the category of \'etale $G$-sets.
$\B(G)$ is equivalent to the presheaf category on $L(G)$.

An ordered groupoid $G$ is said to be {\em principally inductive\/}
if for each identity $e$ the poset
$e^{\downarrow} = \{f \in G_0\colon \: f \leq e \}$
is a meet semilattice under the induced order \cite{L0}.
It is worth noting that if $G$ is an ordered groupoid,
then it is principally inductive
precisely when the left-cancellative category $L(G)$ has pullbacks.
Let $G$ be  a principally inductive groupoid.
Define
\[
C(G) = \{(e,x,f) \in G_0 \times G \times G_0 \colon
  \: \mathbf{d}(x)\leq f, \mathbf{r}(x)\leq e \}
\]
and define a partial binary operation by
$(e,x,f)(f,y,i) = (e, x \circ y, i)$.
Observe that the pseudoproduct $x \circ y$ is defined because
$\mathbf{d}(x), \mathbf{r}(y) \leq f$
and the fact that $G$ is assumed to be principally inductive.
$C(G)$ is a category, and when $G$ is the inductive groupoid of an inverse semigroup, then $C(G)$ is the corresponding Cauchy completion.

An ordered groupoid $G$ is said to be an {\em enlargement\/}
of an ordered groupoid $H$
if $H$ is a full subgroupoid of $G$, an order ideal,
and every object in $G$ is isomorphic to an object in $H$.
Equivalently, $H$ is the full subgroupoid of $G$ spanned
by an open subspace of $G_0$ (in the Alexandrov topology)
intersecting each orbit of $G$ on $G_0$.
This notion is introduced in \cite{L1}.
It is routine to verify that ordered groupoid enlargements of principally inductive groupoids are also principally inductive.
Let $S$ and $T$ be inverse semigroups
with associated inductive groupoids $G(S)$ and $G(T)$.
A {\em bipartite ordered groupoid enlargement\/} of $G(S)$ and $G(T)$
is an ordered groupoid $[G(S),G(T)]$ such that:
it is an enlargement of both $G(S)$ and $G(T)$,
the set of objects of $[G(S),G(T)]$ is the disjoint union of the set
of objects of $G(S)$ and $G(T)$,
and for each $e \in G(S)_0$ there exists an arrow $x$ such that
$\dom(x) = e$ and $\ran(x) \in G(T)_0$, and vice versa.

There is evidently a connection between enlargements and (strong) Morita equivalence since Steinberg \cite{S} observes that if the inverse semigroup $S$
is an enlargement of an inverse semigroup $T$,
then $S$ and $T$ are strongly Morita equivalent,
and Lawson \cite{L1} observes that they are semigroup Morita equivalent.

We shall say that two inverse semigroups, regarded as ordered groupoids,
are {\em enlargement Morita equivalent\/}
if there is an ordered groupoid which is an enlargement of them both.

The main goal of this paper is to prove that these four notions of Morita equivalence are the same.
We shall also study the detailed relationship between the two categories of actions of an inverse semigroup $S$:
the category  $S\text{-}\Se$ of unitary actions and the category {\sl \'Etale\/ }of \'etale actions.
We shall prove in \S~\ref{etale} that
the obvious forgetful functor
\[
U\colon\mbox{\sl \'Etale}\lra S\text{-}\Se\;,\; U(X\ra E)=X\;,
\]
is comonadic.
But more is true: the right adjoint of $U$ is monadic, from which it follows
that $S\text{-}\Se$ is equivalent to $\PSh(C(S))$.
In fact, in \S~\ref{invsemi}
we shall prove that this result generalizes to all semigroups
with right local units,
thus making a direct connection between the Morita equivalence
of semigroups with right local units described in \cite{AL,T1,T2,T3}
and the Morita theory of categories described in \cite{EZ}.\\


\subsubsection*{Acknowledgements}
The authors are grateful to the editor and the referees for their encouraging comments which have led to a better presentation of our ideas.


\section{Morita variants are equivalent}\named{morvar}\setcounter{theorem}{0}

The goal of this section is to prove that the different notions of Morita equivalence that we have defined are in fact the same.
We begin in \S~\ref{ord} by gathering together some basic definitions and
facts about categories that we shall need.

\subsection{Categorical preliminaries}\label{ord}

A {\em weak equivalence\/} from one category to another
is a full and faithful functor that is essentially surjective on objects,
whereas an {\em equivalence\/} is a functor with a pseudo-inverse.
We prefer to distinguish between
weak equivalences and equivalences of categories,
although by the axiom of choice a weak equivalence has a pseudo-inverse.
For instance, an ordered functor $\theta$ that is a local
isomorphism, so that $L(\theta)$ is a weak equivalence (Lemma \ref{Ltheta}),
may not have a pseudo-inverse in the 2-category of ordered groupoids
even though $L(\theta)$ does have one (by choice).
Thus, it is generally good practice to keep track of weak equivalences.
Indeed, in \S~\ref{morvar} we work with weak
equivalences, and ultimately the argumentation does not depend on choice.

We turn to some presheaf preliminaries.
If $\bC$ is a (small) category,
then a contravariant functor from $\bC$ to the category of sets
is called {\em a presheaf}.
Informally, a presheaf is a `right $\bC$-action.'
$\PSh(\bC)$ shall denotes the category of pre\-sheaves on $\bC$.
The functor $Y\colon \bC\lra\PSh(\bC)$ that carries an object $c$
to a representable presheaf $\bC(-,c)$ is full and faithful.
We shall refer to it simply as {\em Yoneda\/} in what follows.
If $P$ is a presheaf on a category $\bC$,
then the {\em category of elements\/} $\bP$ of $P$ is the category
whose objects are pairs $(x,c)$ with $c$ an object of $\bC$ and $x\in P(c)$.
A morphism $f\colon (x,c)\ra (x',c')$ is a morphism $f\colon c\ra c'$
such that $P(f)(x')=x$.
The Yoneda lemma says that an object $(x,c)$ can alternatively be viewed
as a natural transformation $c\arr{x} P$,
where we denote by $c$ the corresponding representable presheaf.
The requirement on $f$ then says that the diagram
\[
\xymatrix{c\ar[rr]^f\ar[rd]_x&&c'\ar[ld]^{x'}\\ &P&}
\]
commutes.
If $\bP\arr{K}\bC$ is the functor sending $(x,c)$ to $c$,
then
\begin{equation}
\colimit {\bP} YK\cong P\;.
\label{YKP}
\end{equation}
(This generalizes the fact that if $M$ is a monoid and $X$ is an $M$-set,
then $X\otimes_M M\cong X$.)
A functor $\bP\arr{F}\bC$ is said to be {\em a discrete fibration\/}
when every morphism $c\arr{m}F(y)$ in $\bC$
has a unique lifting $x\arr{n}y$ to $\bP$.
The isomorphism (\ref{YKP}) is part of the
well-known equivalence between the category of discrete fibrations over $\bC$
and $\PSh(\bC)$~\cite{El}.
The equivalence associates with a presheaf $P$ the discrete fibration $K$
of elements of $P$ just described,
and with a discrete fibration $F$ the colimit $\colimit {\bP} YF$.

We next present some categorical preliminaries
on Morita equivalence of categories.
Details can be found in Chapters 6 and 7 of~\cite{Bo}.
One approach to Morita theory for categories
involves what are called essential points of a topos \cite{Bunge},
whereas another uses what are called profunctors
or bimodules or distributors \cite{Bo}.
It is the second approach we shall use in common with \S~\ref{invsemi}.

Categories $\bA$ and $\bB$ are said to be {\em Morita equivalent\/}
if their presheaf categories are equivalent.
A {\em Morita context\/} for $\bA$ and $\bB$ is a category $\bU$
together with a diagram
$$
\xy
\Vtriangle(0,0)|alr|/`>`>/<300,300>[\bA`\bB`\bU;``]
\endxy
$$
of weak equivalences.

Let $\bC$ and $\bD$ be (small) categories.
A {\em profunctor} $U\colon \bC\pro \bD$ is by definition a functor
\[
U\colon \bD\op\times \bC\ra \Se\;
\]
which can be thought of as a $(\bC,\bD)$-biset.
By exponentiation, this transposes to a functor
$U\colon \bC\lra\PSh(\bD)$,
which in turn corresponds by colimit-extension along Yoneda
to a colimit-preserving functor
\begin{equation}
U\colon \PSh(\bC)\lra\PSh(\bD)\;.
\label{U}
\end{equation}
Categories, profunctors, and natural transformations form a bicategory
(a natural transformation in this context amounts to a biset morphism).
For any $\bC$, the identity profunctor $\bC\ra\bC$ is
the hom-functor $\bC(-,-)$,
which corresponds to Yoneda $\bC\ra \PSh(\bC)$.
Composition of profunctors is given by tensor product.
It is convenient to denote a profunctor $\bC\pro\bD$
and the corresponding functors  $\bD\op\times \bC\ra \Se$,
$\bC\lra\PSh(\bD)$, and (\ref{U}) by one and the same symbol.

We say that a profunctor has a right adjoint if it has
a right adjoint in the usual bicategorical sense.
It follows that a profunctor $\bC\pro \bD$ has a right adjoint
if and only if the corresponding colimit-preserving functor
(\ref{U}) has a colimit-preserving right adjoint (it always has
a right-adjoint, but the right adjoint may not preserve colimits).

Let $\bC$ be a category.
We say that $\bC = [\bA,\bB]$ is
{\em bipartite (with left part $\bA$ and right part $\bB$)}
if it satisfies the following conditions:
\begin{description}
\item[{\rm (B1)}]
$\bC$ has full subcategories
$\bA$ and $\bB$ such that $C_0=A_0\cup B_0$ disjointly.
\item[{\rm (B2)}]
For each object $a\in A_0$
there exists an isomorphism $x$ with domain $a$ and codomain in $B_0$;
for each object $b \in B_0$
there exists an isomorphism $y$ with domain $b$ and codomain in $A_0$.
\end{description}
A bipartite category
$\bC = [\bA,\bB]$ is a disjoint union of four kinds of arrows:
those in $\bA$, those in $\bB$, those starting in $A_0$ and ending in $B_0$,
and those starting in $B_0$ and ending in $A_0$.  Clearly,
$$
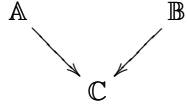
\xy
\Vtriangle(0,0)|alr|/`>`>/<300,300>[\bA`\bB`\bC;``]
\endxy
$$
is a Morita context.

An idempotent $c\arr{e}c$ of a category \emph{splits} if it factors $c\arr{f}r\arr{s}c$, such that $fs=1_r$.
For instance, later we use the fact that idempotents split
in the category $C(S)$ defined in \S~\ref{intro}.

Clearly if two categories have a Morita context,
then they are Morita equivalent.
Our immediate goal is to show that the converse holds
if idempotents split in the two categories,
and moreover, in this case the two categories
have a Morita context coming from a bipartite category.

The following two results are well-known \cite{Bo}.

\begin{lemma}\named{retract}
Suppose that a profunctor $U\colon \bC \pro \bD$ has a right adjoint.
Then for every object $c$ of $\bC$,
$U(c)$ is a retract of a representable in $\PSh(\bD)$.
Moreover, if idempotents split in $\bD$,
then every $U(c)$ is isomorphic to a representable.
\end{lemma}

A presheaf is said to be {\em indecomposable\/}
if the covariant hom-functor associated with it
preserves coproducts.

\begin{proposition}\named{projconn}
A presheaf on a small category $\bC$ is indecomposable and projective iff
it is a retract of a representable.
If idempotents split in $\bC$,
then a presheaf is indecomposable and projective
if and only if it is isomorphic to a representable.
\end{proposition}


An {\em equivalence profunctor} is a profunctor that is an
equivalence in the bicategory of profunctors.
In algebraic terms, an equivalence amounts to a $(\bC,\bD)$-biset $U$
and a $(\bD,\bC)$-biset $V$ such that
\[
U\otimes_{\bD} V\iso \bC(-,-)\qquad V\otimes_{\bC} U\iso \bD(-,-)\;.
\]
It is known~\cite{Bo} that $\PSh(\bC)$ is equivalent to $\PSh(\bD)$
if and only if there is an equivalence profunctor $U\colon \bC\ra \bD$.
Indeed, $U\colon \PSh(\bC)\lra \PSh(\bD)$ is an equivalence of categories
if and only if the corresponding profunctor is
an equivalence profunctor~\cite{Bo}.

We sometimes denote the coproduct of two sets $A$ and $B$ by $A+B$,
commonly understood as `disjoint union.'

\begin{proposition}\named{CED}
Suppose that idempotents split in both $\bC$ and $\bD$.
An equivalence $U\colon \PSh(\bC)\lra \PSh(\bD)$,
i.e., an equivalence profunctor $U\colon \bC\ra \bD$,
gives rise to a Morita context
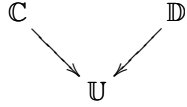
\[
\xy
\Vtriangle(0,0)|alr|/`>`>/<300,300>[\bC`\bD`\bU;``]
\endxy
\]
such that $\bU=[\bC,\bD]$.
\end{proposition}
\proof
We define a category $\bU$ as follows.
Let $U_0=C_0+D_0$,
and let $U_1=C_1+D_1+X$,
where $X$ is the collection of all
natural transformations between objects
$U(c)$ and $d$ in $\PSh(\bD)$
(as usual, we omit notation for both Yoneda functors).
For instance, a natural transformation $U(c)\ra d$ is a morphism
$c\ra d$ in $\bU$.
Then $\bU$ is a category,
and by Lemma \ref{retract} we have $\bU=[\bC,\bD]$.
\qed

\subsection{Topos equivalence implies strong equivalence}\label{tese}


Let $S$ and $T$ be inverse semigroups,
and assume that the toposes $\B(S)$ and $\B(T)$ are equivalent.
We use Proposition\ \ref{CED} to show that $S$ and $T$
are strongly Morita equivalent.
In this case, $\bC=L(S)$ and $\bD=L(T)$ are left-cancellative
categories, so the identities are their only (split) idempotents.
By Proposition\ \ref{CED} (and its proof), there is an equivalence
$U\colon \B(S)\equ\B(T)$ if and only if there is a Morita context
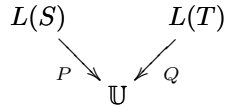
\[
\xy
\Vtriangle(0,0)|alr|/`>`>/<300,300>[L(S)`L(T)`\bU;`P`Q]
\endxy
\]
where $\bU$ is the (left-cancellative)
category whose objects are the idempotents
of $S$ and $T$ (disjoint collection).
$\bU=[L(S),L(T)]$ has three kinds of morphisms:
\begin{enumerate}
\item
those of $L(S)$,
\item
those of $L(T)$, and
\item
the connecting ones between $d\in E(S)$ and $e\in E(T)$,
which are understood as natural transformations
between presheaves $U(d)$ and $Y(e)$ in $\B(T)$,
where $U\colon L(S)\pro L(T)$ is the equivalence profunctor and $Y$ is Yoneda.
\end{enumerate}

We may reorganize this data into an equivalence biset in the semigroup sense.
In what follows, we do not distinguish notationally
between the object $e$ of $L(T)$ and the presheaf $Y(e)$.
Let $X$ denote the set of connecting isomorphisms
from an idempotent of $T$ to an idempotent of $S$;
that is, the morphisms of type 3 above,
but only the isomorphisms and only in the direction from $T$ to $S$.

The action by $S$ is precomposition,
which we write as a left action.
Let $e\arr{x}d$ be an element of $X$: this is
an isomorphism $x\colon e\iso U(d)$ in $\B(T)$.
Let $s\in S$. If $s^*s=d$, then $sx$ is
the composite isomorphism $e\iso U(d)\iso U(ss^*)$,
i.e., $U(ss^*,s)x$.
This defines a partial action by $S$,
which we can make total with the help of the following lemma.

\begin{lemma}\named{ba}
Let $U\colon \B(G)\equ \B(H)$ be an equivalence of classifying toposes
of ordered groupoids $G$ and $H$.
Let $b\leq d$ in $G_0$ and $x\colon e\iso U(d)$ be an
isomorphism of $\B(H)$.
Then there is a unique idempotent $a\leq e$ in $H_0$,
and a unique isomorphism $bx\colon a\iso U(b)$ such that
\[
\xy
\square(0,0)|alra|/>`>`>`>/<400,300>[a`U(b)`e`U(d);bx```x]
\endxy
\]
is a pullback in $\B(H)$.
\end{lemma}
\proof
By Lemma \ref{retract},
there is $c\in H_0$ and an isomorphism $y\colon c\iso U(b)$.
Consider the composite
$$
c\iso U(b)\ra U(d)\iso e
$$
in $\B(H)$, where the last isomorphism is $x^{-1}$.
By Yoneda,
this comes from a unique morphism $c\arr{t} e$ in $L(H)$.
Let $a=\mathbf{r}(t)\leq e$, and $bx=yt^{-1}$.

Such an $a$ is unique because a subobject
(which is an isomorphism class of monomorphisms)
of a representable $e$ corresponds uniquely to a
downclosed subset of elements of $H_0$ under $e$,
and a principal one corresponds uniquely to
an element of $H_0$ under $e$.
If $a$ and $a'$ both make the square a pullback,
then they are in the same isomorphism class of
monomorphisms into $e$, hence they represent the same subobject, hence $a=a'$.
The isomorphism $bx$ is also unique
because $U(b)\ra U(d)\iso e$ is a monomorphism.
\qed
Returning to inverse semigroups,
we see how to make the action total:
let $b=ds^*s\leq d$, and let $sx=sd\cdot bx$.

The inner product $\<\;,\;\>\colon X\times X\ra S$ is defined as follows.
If two isomorphisms $x\colon e\iso U(d)$ and $y\colon e\iso U(c)$ have
the same domain, then $\<x,y\>=yx^{-1}$.
This is an isomorphism of $\B(T)$ between $U(d)$ and $U(c)$,
but that amounts to an isomorphism of $L(S)$,
which in turn is precisely an element of $S$.
In general, the inner product of $x\colon f\iso U(d)$ and $y\colon e\iso U(c)$
is defined by using variations of Lemma \ref{ba}.
$$
\xy
\square(0,0)|alra|/>`>`>`>/<400,300>[U(a)`ef`U(d)`f;```x^{-1}]
\square(900,0)|alra|/>`>`>`>/<400,300>[ef`U(b)`e`U(c);```y]
\endxy
$$
These ``variations'' can be established in the same
way as in Lemma \ref{ba}, or they can be deduced from
Lemma \ref{ba} by transposing under the pseudo-inverse $V$ of $U$.
For example, the right-hand square above can be obtained
by applying Lemma \ref{ba} (with $V$ instead of $U$)
to the transpose of $y^{-1}$, as in the following diagram.
$$
\xy
\square(0,0)|alra|/>`>`>`>/<400,300>[b`V(ef)`c`V(e);```\widehat{y^{-1}}]
\endxy
$$

The right action by $T$ and the inner product
$[\;,\;]\colon X\times X\ra T$ are entirely analogous.
The axioms (M1) - (M7) are easily verified.
For example, for any $x\colon f\iso U(d)$,
the rule (M3) $\<x,x\>x=x$ is the fact that the composite $xx^{-1}x$
is equal to $x$ (in $\bU$):
$$
f\iso U(d)\iso f\iso U(d)\;\;;\; \<x,x\>x=xx^{-1}x=x\;.
$$

\subsection{Strong equivalence implies topos equivalence}\label{sete}

Although Steinberg \cite{S} proves this (assuming choice),
it may be of interest to see how to build a Morita context
$$
\xy
\Vtriangle(0,0)|alr|/`>`>/<300,300>[L(S)`L(T)`\bU;`P`Q]
\endxy
$$
in the category sense from an equivalence biset $X$.

By definition, the objects of the bipartite category $\bU=[L(S),L(T)]$ are disjointly the objects of $L(S)$ and $L(T)$,
which are the idempotents of $S$ and of $T$.
A morphism of $\bU$ is either:
\begin{enumerate}
\item
one of $L(S)$,
\item
one of $L(T)$,
\item
one of the form $(x,d)\in X\times E(S)$,
such that $\<x,x\>\leq d$,
where the domain of this morphism is $[x,x]\in E(T)$,
and its codomain is $d$, or
\item
one of the form $(x,e)\in X\times E(T)$,
such that $[x,x]\leq e$,
where the domain of this morphism is $\<x,x\>\in E(S)$,
and its codomain is $e$.
\end{enumerate}
We compose the various kinds of morphisms in $\bU$
by using the inner products and actions in $X$ by $S$ and $T$.
For example, by definition
\[
\xy
\qtriangle(0,0)|alr|/>`>`>/<300,300>[s^*s`d`e;s`s^*x`x]
\endxy
\]
commutes in $\bU$,
where $s\in S$, $d\in E(S)$, $x\in X$, $d=\<x,x\>$,
$s=ds$, $e\in E(T)$ and $[x,x]\leq e$.
In other words, we define $(x,e)(s,d)=(s^*x,e)$.
The pair $(s^*x,e)$ is indeed a legitimate morphism of $\bU$ because
the idempotent product $[x,x][s^*x,s^*x]$ is equal to
\[
[x,\<x,s^*x\>s^*x]=[x,\<x,x\>ss^*x]=[x,dss^*x]
=[x,ss^*x]=[s^*x,s^*x]\;.
\]
Therefore, $[s^*x,s^*x]\leq[x,x]\leq e$.
The domain of $(s^*x,e)$ is
\[
\<s^*x,s^*x\>=s^*\<x,x\>s=s^*ds=s^*s\;,
\]
which is the domain of $(s,d)$ as it should be.
For another example,
\[
\xy
\qtriangle(0,0)|alr|/>`>`>/<400,300>[\<x,x\>`{[y,y]}`e;x`\<y,x\>`y]
\endxy
\]
commutes, where $[x,x]\leq [y,y]$.
The domain of the composite $\<y,x\>$ is
\[
\<y,x\>^*\<y,x\>=\<x,y\>\<y,x\>=\<x{[y,y]},x\>=\<x,x\>\;,
\]
since $x=x[x,x]=x[x,x][y,y]=x[y,y]$.
It follows that
$\bU$ is a category, that $\bU=[L(S),L(T)]$,
and that the obvious functors $P,Q$ are weak equivalences.

\begin{corollary}\named{Ulc}
The category $\bU$ constructed from an equivalence biset
is left-cancellative.
\end{corollary}
\proof
This is true because $\bU$ is weakly equivalent to a left-cancellative
category.
However, the following calculations give more information.
For example, if
\[
\xy
\qtriangle(0,0)|alr|/>`>`>/<300,300>[s^*s`d`e;s`y`x]
\endxy
\]
commutes in $\bU$, where $d=\<x,x\>$ and $[x,x]\leq e$,
then $y=s^*x$ (by definition) and
\[
s=ds=\<x,x\>s=\<x,s^*x\>=\<x,y\>\;.
\]
Thus, $s$ is uniquely determined by $x$ and $y$.
The other possibility, but keeping $(x,e)$, is
\[
\xy
\qtriangle(0,0)|alr|/>`>`>/<400,300>[{[y,y]}`d`e;y`t={[x,y]}`x]
\endxy
\]
where $\<y,y\>\leq d$.
Then $y$ is determined by $x$ and $t$ since
\[
y=\<y,y\>y=\<x,x\>\<y,y\>y=\<x,x\>y=x{[x,y]}=xt\;.
\]
It follows that $(x,e)$ is a monomorphism.\qed


\subsection{Topos equivalence and enlargement equivalence}

In this section, it is no harder to work with ordered groupoids
more general than inductive groupoids.

A poset map $P\arr{f}Q$ is said to be a discrete fibration (\S~\ref{ord})
if for every $x\leq f(y)$ in $Q$ there is a unique $z\in P$ such that $f(z)=x$.
For example, the domain map of an ordered groupoid
is by definition a discrete fibration.
A poset map is a discrete fibration if and only if it is \'etale
(i.e., a local homeomorphism) for the Alexandrov topology.

An ordered functor $\theta:G \ra H$ is said to be
{\em a local isomorphism\/}
if it satisfies the following two conditions.

\begin{description}
\item[{\rm (LI1)}]
the underlying groupoid functor of $\theta$ is a weak equivalence;
\item[{\rm (LI2)}]
the object function
$\theta_0\colon G_0\ra H_0$ is a discrete fibration of posets.
\end{description}
An enlargement is a local isomorphism.

\begin{lemma}\named{Ltheta}
An ordered functor $\theta\colon G\ra H$ is a local isomorphism
if and only if
$L(\theta)\colon L(G)\ra L(H)$ is a weak equivalence.
\end{lemma}
\proof
Assume $\theta$ is a local isomorphism.
Clearly $L(\theta)$ is essentially surjective if $\theta$ is.
$L(\theta)$ is full: let $\theta(d)\arr{t}\theta(e)$ be a morphism of $L(H)$.
Consider the unique lifting $c\leq e$ of $\mathbf{r}(t)\leq \theta(e)$,
so that $\theta(c)=\mathbf{r}(t)$.
Since $\theta$ is full there is $d\arr{s}e$ (in $G$) such that
$\theta(s)=t$.
Thus, $L(\theta)(s)=t$.
$L(\theta)$ is faithful: suppose that $L(\theta)(s)=L(\theta)(t)$,
where $s,t\colon  d\ra e$ in $L(G)$.
Let $c=\theta(\mathbf{r}(s))=\theta(\mathbf{r}(t))$.
The two inequalities $\mathbf{r}(s)\leq e$ and $\mathbf{r}(t)\leq e$
both lie above $c\leq \theta(e)$,
so they must be equal by the uniqueness of liftings along $\theta_0$.
Thus, if $\theta$ is faithful, then $s=t$.

For the converse, if $L(\theta)$ is a weak equivalence,
then we see easily that $\theta$ satisfies (LI1).
One can verify (LI2) directly,
but we prefer the following more conceptual argument.
We have a commuting square of toposes
$$
\xy
\square(0,0)|alra|/>`>`>`>/<600,300>[\PSh(G_0)`\PSh(H_0)`\B(G)`\B(H);```]
\endxy
$$
where the bottom horizontal is the equivalence associated with
the weak equivalence $L(\theta)$.
Since the two geometric morphisms depicted
vertically are \'{e}tale, so is the top horizontal.
Therefore, $G_0\ra H_0$ is a discrete fibration.
\qed

\begin{theorem}\named{united}
The following are equivalent for ordered groupoids $G$ and $H$:
\begin{enumerate}
\item
the classifying toposes of $G$ and $H$ are equivalent;
\item
$G$ and $H$ have a joint bipartite enlargement $[G,H]$;
\item
there is an ordered groupoid $K$ and local isomorphisms
$G\ra K\la H$.
\end{enumerate}
\end{theorem}
\proof
$(1)\2cell (2)$.  Given an equivalence $U\colon \B(G)\equ \B(H)$,
consider the groupoid $K$
such that $K_0=G_0+H_0$ and $K_1=G_1+H_1+Y$,
where $Y$ is set of isomorphisms of $\B(H)$ between objects
$U(d)$ and $e$.
$K_1$ is partially ordered:  for $i\colon U(d)\iso e$ and $j\colon U(a)\iso b$,
we declare $i\leq j$ when $d\leq a$ in $G_0$ and $e\leq b$ in $H_0$
and the square of natural transformations
$$
\xy
\square(0,0)|alra|/>`>`>`>/<400,300>[U(d)`e`U(a)`b;i```j]
\endxy
$$
commutes in $\B(H)$.
The definition of $\leq$
for isomorphisms in the other direction is similar.
By Lemma \ref{ba}, the domain map $K_1\ra K_0$ is a discrete fibration.

$(2)\2cell (3)$ holds because an enlargement is a local isomorphism.

$(3)\2cell (1)$ holds because
given such local isomorphisms, then $\B(G)$ and $\B(H)$ are
equivalent by Lemma \ref{Ltheta} since the geometric morphism
associated with a weak equivalence of categories is an equivalence.
\qed

We construct from a given equivalence biset $X$ between inverse semigroups
$S$ and $T$ a common ordered groupoid enlargement of
$G(S)$ and $G(T)$, denoted $G(S,T;X)\,$.
We do this again in Theorem~\ref{boge} using semigroup methods.
We start with the presheaf
$$
\bfS(e) = \left\{
\begin{array}{lr}
\{s\in S\mid s^*s=e\}+\{x\in X\mid \<x,x\>=e\}\;, & e\in E(S)\\
\{t\in T\mid t^*t=e\}+\{x\in X\mid [x,x]=e\}\;, & e\in E(T)
\end{array}\right.
$$
on the left-cancellative category $\bU$ built from $X$ (as in Cor.\ \ref{Ulc}).
Let $\bS_0\ra \bU$ denote the discrete fibration
corresponding to the presheaf $\bfS$.
$\bS_0$ is the category of elements of $\bfS$,
whose objects are `elements' $e\arr{u}\bfS$.
The category of elements of any presheaf
on a left-cancellative category is a preorder,
so that $\bS_0$ is a preorder.
The category pullback
$$
\xy
\square(0,0)|alra|/>`>`>`>/<500,300>[\bS_1`\bS_0`\bS_0`\bU;```]
\endxy
$$
defines a {\em pre\/}ordered groupoid $(\bS_0,\bS_1)$.
Let $G(S,T;X)$ denote the posetal collapse of $(\bS_0,\bS_1)$:
the object-poset of $G(S,T;X)$ equals the posetal collapse of $\bS_0$,
which may be identified with the map
$$
\bS_0\epi E(S)+E(T)
$$
such that an element
$$
e\arr{u}\bfS\mto
\left\{
\begin{array}{ll}
uu^* & u\in S\;{\rm or}\; u\in T\\
\<u,u\> & u\in X\;{\rm and}\; e=[u,u]\\
{[u,u]} & u\in X\;{\rm and}\; e=\<u,u\>
\end{array}\right.\;.
$$
Likewise, the morphism-poset of $G(S,T;X)$
equals the posetal collapse of $\bS_1$.
Moreover, the underlying groupoid of $G(S,T;X)$,
where we ignore its order structure,
equals the isomorphism subcategory of $\bU$.

To conclude this section, we shall relate the strong Morita equivalence of two inverse semigroups
with the two categories $L(S)$ and $C(S)$ that we have defined for any inverse semigroup $S$.

\begin{lemma}\named{liwe}
Let $G$ and $H$ be principally inductive.
Then an ordered functor
$\theta\colon G\ra H$ is a local isomorphism
if and only if $C(\theta)\colon C(G)\ra C(H)$ is a weak equivalence.
\end{lemma}
\proof
The forward implication is similar to the proof of Lemma \ref{Ltheta}.
On the other hand, if $C(\theta)$ is a weak equivalence,
then so is $L(\theta)$
because $L(G)$ equals the subcategory
of $C(G)$ consisting of those morphisms with retracts~\cite{F}.
We may now appeal to Lemma \ref{Ltheta}.
\qed

\begin{proposition}\named{CGCH}
Let $G$ and $H$ be principally inductive ordered groupoids.
Then the following are equivalent:
\begin{enumerate}
\item the classifying toposes of $G$ and $H$ are equivalent;
\item
the categories $L(G)$ and $L(H)$ form a Morita context.
\item
the categories $C(G)$ and $C(H)$ form a Morita context;
\end{enumerate}
\end{proposition}
\proof
(1) and (2) are equivalent because idempotents split in
the left-can\-cellative category $L(G)$,
and since $\B(G)\equ\PSh(L(G))$.

(2) and (3) are equivalent because
$C(G)$ is canonically equivalent to the category ${\rm Span}(L(G))$,
where the Span of a category with pullbacks is given by
the same objects, but whose morphisms
are spans $\cdot \la \cdot \ra\cdot$ in the given category.
Hence, a Morita context for $C(G)$ and $C(H)$
comes from one for $L(G)$ and $L(H)$ by applying the Span construction.
(This aspect is further explained following Lemma \ref{projind}.)
Conversely,
a Morita context for $L(G)$ and $L(H)$ can be obtained
from one for $C(G)$ and $C(H)$
because as in the proof of Lemma \ref{liwe} $L(G)$ equals the
retract subcategory of $C(G)$.
\qed

\subsection{Strong equivalence and semigroup equivalence}\named{invsemi}

We shall prove that strong Morita equivalence and semigroup equivalence are the same.
But to do this we shall prove a theorem for a much wider class of semigroups than just the inverse ones.
We recall that if $S$ is a semigroup with right local units,
then  $S\text{-}\Se$ denotes the category of closed right $S$-sets.

\begin{lemma} Let $S$ be a semigroup with right local units.
Then the category  $S\text{-}\Se$ has all small colimits,
and they are created by the underlying set functor.
\end{lemma}
\proof
Let $\Se_S$ be the category of sets with a right action by $S$.
It is well-known that $\Se_S$ is complete and cocomplete,
and that limits and colimits are created by the underlying set functor.
The functor $\Se_S\lra \Se_S$ given by $X\mapsto X\otimes_S S$
(with the usual action) has a right adjoint $X\mapsto \hom_S(S,X)$,
so that it therefore preserves colimits.
The collection of morphisms $\mu_X\colon X\otimes_S S\ra X$
given by $x\otimes s\mapsto xs$ constitute a natural transformation
$\mu$ from $(-)\otimes_S S$ to the identity functor on $\Se_S$,
and $S\text{-}\Se$ is the full subcategory of $\Se_S$ on the objects for which $\mu$ is an isomorphism.
It follows that $S\text{-}\Se$ is closed under small colimits.
Indeed, if $\bD$ is a small category and $F\colon \bD\ra S\text{-}\Se$
is a functor, then writing $F\otimes _S S$ for the functor
$d\mapsto F(d)\otimes_S S$,
we have that $F\otimes_S S\iso F$ as functors to $\Se_S$ via the natural transformation with components $\mu_{F(d)}$.
Thus
\[
\colimit{\bD} F\iso \colimit{\bD} (F\otimes_S S)\iso
 (\colimit \bD F)\otimes_S S
\]
since tensor product commutes with colimits.
Diagram chasing reveals that the isomorphism is given by $\mu$.
\qed

As usual, $Y$ denotes the Yoneda functor $C(S)\lra\PSh(C(S))$.
There is also a functor $F:C(S)\lra S\text{-}\Se$
defined as follows: for each idempotent $e$ in $S$,
corresponding to the identity $(e,e,e)$, we define $F(e) = eS$,
and if $(f,a,e)$ is an arrow in $C(S)$ from $e$ to $f$,
then $F(f,a,e) \colon eS \ra fS$ is given by $x \mapsto ax$.
This is a well-defined functor because $eS$ really is a closed right $S$-set.
The proof of this follows by an argument similar to that used in \cite{AL}.

\begin{theorem}\label{directproof}
Let $S$ be a semigroup with right local units.
Then the categories $S\text{-}\Se$ and $\PSh(C(S))$ are equivalent.
\end{theorem}
\proof Let $S$ be a semigroup with right local units.
We may easily define a functor $Q$ from  $S\text{-}\Se$
to $\PSh(C(S))$ as follows.
If $X$ is a closed right $S$-set, then $Q(X)$ is the presheaf on $C(S)$
defined by $Q(X)(e)=Xe\,$.
The transition map of $Q(X)$ for a morphism $(e,s,f)$ of $C(S)$
is given by $Q(X)(e,s,f)(x)=xs$, which we more conveniently denote by $x(e,s,f)$.
The restriction of an $S$-equivariant map $X\arr{h}Y$ to $e$ gives the component at $e$ of a natural transformation $Q(X)\arr{Q(h)}Q(Y)$.
The following diagram commutes.
\[
\xy
\Atriangle(0,0)|lrb|/>`>`>/<400,400>%
  [C(S)`S\text{-}\Se`\PSh(C(S));F `Y`Q]
\endxy
\]
We claim that $Q$ has a left adjoint $Q_!$,
which is defined by the colimit extension:
\[
Q_!(P)=\colimit{}\left( \bP\ra C(S)
  \arrow{F}S\text{-}\Se\right)\;,
\]
where $\bP\ra C(S)$ is the discrete fibration corresponding to a presheaf $P$.

To show that the adjunction $Q_!\adj Q$ is an (adjoint) equivalence,
it suffices to show that $Q$ is full, faithful,
and that for any presheaf $P$, the unit $P\ra Q(Q_!(P))$ is an isomorphism.

\begin{Claim}  $Q$ preserves small colimits.\end{Claim}
\proof
$Q$ clearly preserves coproducts since they set-theoretic
in $S\text{-}\Se$ and componentwise in $\PSh(C(S))$.
$Q$ also preserves coequalizers.
The coequalizer of two morphisms
\[
\xy
\morphism(0,0)|a|/{@{>}@/^1em/}/<500,0>[X`Y;f]
\morphism(0,0)|b|/{@{>}@/_1em/}/<500,0>[\phantom{X}`\phantom{X};g]
\endxy
\]
in $S\text{-}\Se$ is created by the underlying set functor and hence is the set $Y/R$,
where $R$ is the equivalence relation
generated by identifying $f(x)$ with $g(x)$ for $x\in X$.
This is preserved by $Q$ since
if $ye=y'e$ and $y=y_1,\ldots,y_m=y'$ is a zig-zag of elements,
so that for each $i$ there is an $x_i\in X$
such that either $f(x_i)=y_i$ and $g(x_i)=y_{i+1}$ or vice versa,
then $y=ye=y_1e,\ldots,y_me=y'e=y'$ is a zig-zag,
which proves that $x, y$ get identified in the quotient of $Ye$
obtained when constructing the coequalizer of $Q(f),Q(g)$
in $\PSh(C(S))$.
Conversely, an identification in $Ye$ when forming
the coequalizer of $Q(f)$ and $Q(g)$
yields an identification of the corresponding elements in $Y$.
\qed

\begin{Claim}
$Q$ is faithful.
\end{Claim}
\proof
If $f,g\colon  X\ra Y$ are two morphisms with $Q(f)=Q(g)$,
then for any idempotent $e$, $f$ and $g$ agree on $Xe$.
But $X$ is the union of the $Xe$ over all $e$, so $f=g$.
Thus $Q$ is faithful.
\qed

Our next claim is where we use that the action is closed.

\begin{Claim}
$Q$ is full.
\end{Claim}
\proof
Let $Q(X)\arr{h}Q(Y)$ be a natural transformation.
Then we define a map $H\colon X\times S\ra Y$
by $H(x,s)=h_e(xs)$, where $e$ is any idempotent such that $se=s$.
This is well-defined because if $se'=s$ and $f\in E(S)$ satisfies $xf=x$,
then $h_e(xs) = h_e(x(f,fs,e)) = h_f(x)(f,fs,e)=
h_f(x)s=h_f(x)(f,fs,e') = h_{e'}(x(f,fs,e'))=h_{e'}(xs)$.

Next observe that $H$ satisfies $H(xs,t)=H(x,st)$
for all $x\in X$ and $s,t\in S$.
Indeed, if $t=te$ with $e\in E(S)$,
then $st=ste$ so that $H(x,st)= h_e(xst) = H(xs,t)$.
Thus there is a well-defined induced map
$H\colon X\otimes_S S\ra Y$ given by $x\otimes s\mapsto h_e(xs)$,
where $se=s$ with $e\in E(S)$.
Observe that $H$ is an $S$-set morphism because
if $se=s,tf=t$ with $e,f\in E(S)$,
then  $H(x\otimes s)t=h_e(xs)t = h_e(xs)(e,et,f)
=h_f(xs(e,et,f)) = h_f(xst) = H(x\otimes st) = H((x\otimes s)t)$.

Let $H'\colon X\ra Y$ be the composition $H\mu^{-1}$,
where $\mu\colon X\otimes S\ra X$ is the canonical isomorphism.
Then for $x\in Xe$, we have
\[
Q(H')_e(x) = H'(x)=H(x\otimes e)=h_e(x)\;,
\]
so that $Q(H')=h$ establishing that $Q$ is full.
\qed

Finally, we show that the unit for $Q_!\adj Q$ is an isomorphism.
Let $P$ be a presheaf on $C(S)$ with corresponding
category of elements $\bP\arr{K}C(S)$.
We have
\[
\colimit{\bP} Y\cdot K\iso P\;,
\]
where $Y$ denotes the Yoneda functor.
Since $Q$ preserves small colimits,
we have
\[
P\iso \colimit{\bP}Y\cdot K\iso
\colimit{\bP}Q\cdot F\cdot K \iso
Q(\colimit{\bP} F\cdot K)\iso Q(Q_!(P))\;.
\]
This isomorphism is the unit $P\ra Q(Q_!(P))$.
\qed

As a corollary we obtain the analogue of a result
proved by Lawson for Morita equivalence of semigroups
with local units~\cite{AL},
which is again analogous to the results for monoids and categories.

\begin{corollary}
If $S$ and $T$ are semigroups with right local units,
then $S$ and $T$ are Morita
equivalent if and only if there is a Morita context for $C(S)$ and $C(T)$.
\end{corollary}
\proof
This follows from the Theorem~\ref{directproof}
since $C(S)$ and $C(T)$ have split idempotents.
\qed

Talwar~\cite{T2} considers a more general notion of a closed $S$-set
for semigroups satisfying $S^2=S$.
Here an $S$-set $X$ is {\em closed\/} if the natural
morphism $\hom_S(S,X)S\otimes S\ra S$ given by $\alpha t\otimes s=\alpha(ts)$
is an isomorphism, where $\hom_S(S,X)$ is the set of $S$-equivariant
maps from $S$ to $X$.
Denote the corresponding category by $S\text{-}\Se\/$.
If $S$ has local units, he shows that this is equivalent to
the previous notion of closed $S$-set~\cite{T1}.
Talwar calls $S$ a {\em sandwich semigroup\/} if $S=SE(S)S$,
and he proves that $S\text{-}\Se\/$ is equivalent
to $T\text{-}\Se\/$~\cite{T2}, where $T=E(S)SE(S)$.
Of course $T$ has local units.  Also $C(S)=C(T)$.
If $S$ is finite, then $S=S^2$ if and only if $S=SE(S)S$.
Our results have the following corollary.

\begin{corollary}
Let $S$ be a sandwich semigroup.
Then $S\text{-}\Se\/$ is equivalent to $\PSh(C(S))$.
Consequently, if $S$ and $T$ are sandwich semigroups,
then $S\text{-}\Se\/$ is equivalent to $T\text{-}\Se\/$
if and only if there is a Morita context for $C(S)$ and $C(T)$.
\end{corollary}

Finally, we may conclude our proof of the equivalence
between the four types of Morita equivalence defined in \S~\ref{intro}.

\begin{theorem}
Let $S$ and $T$ be inverse semigroups.
Then $S$ and $T$ are strongly Morita equivalent
if and only if they are semigroup Morita equivalent.
\end{theorem}
\proof
In \S~\ref{tese} and \S~\ref{sete} we proved that strong Morita equivalence
is the same as topos Morita equivalence.
In Proposition~\ref{CGCH}, we proved that $S$ and $T$ are topos Morita equivalent
if and only if $C(S)$ and $C(T)$ form a Morita context.
Since the idempotents of $C(S)$ and $C(T)$ split,
they form a Morita context if and only if $\PSh(C(S))$ is equivalent to $\PSh(C(T))$~\cite[Theorem~7.9.4]{Bo}.
Theorem~\ref{directproof} implies $\PSh(C(S))\equ \PSh(C(T))$
if and only if $S$ and $T$ are semigroup Morita equivalent.
\qed

\section{Unitary actions and \'{e}tale actions}\named{etale}

Our goal in this section is to describe in detail the connection between the categories  $S\text{-}\Se\/$ and  {\sl \'Etale\/} in the inverse case.
We have already seen that $S\text{-}\Se\/$ is equivalent
to the presheaf topos $\PSh(C(S))$ (Thm.\ \ref{directproof});
however, it may be illuminating to revisit this fact and several other related ones
in terms of the connection between $S\text{-}\Se\/$ and {\sl \'Etale\/},
without appealing to Thm.\ \ref{directproof}.

\begin{lemma}\named{Scolimits}
$S\text{-}\Se\/$ has all small colimits, created in the category of sets.
All small limits also exist in $S\text{-}\Se\/$
{\em(but they are not created in sets)}.
\end{lemma}
\proof
A small coproduct $\coprod_A X_a$ of unitary actions
is an $S$-set in the obvious way,
which is easily seen to be unitary.
The set coequalizer
\[
\xy
\morphism(500,0)|a|/>>/<500,0>[\phantom{X}`Z;]
\morphism(0,0)|a|/{@{>}@/^1em/}/<500,0>[X`Y;]
\morphism(0,0)|b|/{@{>}@/_1em/}/<500,0>[\phantom{X}`\phantom{X};]
\endxy
\]
of two $S$-maps also has an action by $S$ in an obvious way
(use the universal property of $Z$), which again is unitary.

Limits are slightly more complicated than colimits.
For example, a product $X\times Y$ has underlying set
$\{(x,y)\mid \exists e\in E,  ex=x, ey=y\}$.
Arbitrary products follow the same pattern.
Equalizers, like coequalizers, are created in sets.
\qed

An $S$-set is said to be
{\em indecomposable\/} if its covariant hom-functor
preserves coproducts, or equivalently it cannot be
expressed as a coproduct of two proper sub-$S$-sets.

\begin{lemma}\named{inpro}
An $S$-set $eS$ with $e\in E(S)$ is unitary.
A unitary $S$-set is indecomposable and projective
if and only if it is isomorphic to $eS$,
for some idempotent $e$.
The usual functor
\[
F:C(S)\lra S\text{-}\Se\;,\; F(e)=eS\;,
\]
is full and faithful,
giving a weak equivalence of $C(S)$
with the full subcategory of $S\text{-}\Se\/$
on the indecomposable projectives.
\end{lemma}
\proof
We have seen in Lemma \ref{Scolimits}
that $S\text{-}\Se\/$ has arbitrary coproducts, which are created in \Se.
It can be proved, using essentially the same argument as that in \cite{B},
that in this category epimorphisms are precisely the surjections.
An $S$-set $eS$ is clearly unitary,
and it can be directly verified that it is an indecomposable projective.
Indeed, there is a well-known isomorphism of functors $S\text{-}\Se(eS,-)$
and $(-)e$ given on an $S$-set $X$ by
\begin{equation}\label{natisoproj}
\begin{split}
\eta_X\colon S\text{-}\Se(eS,X)\ra Xe\\ f\mto f(e).
\end{split}
\end{equation}
Note that $\eta_X^{-1}(x)$ is the map $eS\ra X$ given by $s\mto xs$.
Trivially, the functor $X\mapsto Xe$ preserves finite coproducts
and surjective morphisms, whence $eS$ is an indecomposable projective.

Let $X$ be an arbitrary unitary $S$-set, and let $x \in X$.
By unitary,
there exists $s \in S$ and $y \in X$ such that $ys = x$.
Then $xs^*s = yss^*s=ys=x$.
Let \[R(X)=\{(x,e)\in X\times E\mid xe=x\}.\]
The coproduct $\coprod_{R(X)} eS$
is projective and unitary.  Given $(x,e)\in R(X)$, there is a morphism $\pi_{(x,e)}\colon eS\ra X$ with $\pi_{(x,e)}(e)=x$, namely put $\pi_{(x,e)}=\eta_X^{-1}(x)$.  The map $\pi\colon \coprod_{R(X)}eS\ra X$,
given on the component $(x,e)$ by $\pi_{(x,e)}$,
is then a surjection because if $xe=x$, then $\pi_{(x,e)}(e)=x$.

By the same argument as in Proposition~II.14.2 of \cite{M},
every surjection onto a projective is a retraction.
Let $X$ be an arbitrary indecomposable projective.
Thus, the surjection $\pi$ above is a retraction,
so that there is an injective $S$-equivariant map
$\sigma\colon X \ra \coprod_{R(X)} eS$
such that $\pi\cdot\sigma$ is the identity on $X$.
Since $X$ is indecomposable,
$\sigma\colon X \ra eS$ for some $(x,e)$ such that $xe=x$.
This map is necessarily injective.
Since $\pi\cdot\sigma = 1_{X}$
we find that $X =\pi(eS)$, so that $X$ is a cyclic $S$-set.
Therefore, $X\iso \sigma(X)$ and $\sigma(X)$
is a cyclic sub-$S$-set of $eS$, whence a principal right ideal of $S$.
But principal right ideals in an inverse semigroup are generated by idempotents.
Finally, the functor $F(e)=eS$ is full and faithful
because $S\text{-}\Se(dS,eS)\iso eSd=C(S)(d,e)$ by \eqref{natisoproj}.
\qed

We now turn to the category {\sl \'Etale\/}.
Recall that an object of this category
is a set $X$ equipped with a right action by $S$
and a map $X\arr{p}E$ (the \'etale structure)
such that $p(xs)=s^*p(x)s$ and $xp(x)=x$.
Maps in {\sl \'Etale\/} commute with the actions and with the projections to $E$.
Thus, {\sl \'Etale\/} is the full subcategory of $S$-\Se$/E$ on
those objects $X\arr{p}E$ satisfying $xp(x)=x$,
whose inclusion has a right adjoint denoted $V$ in (\ref{VER}).

Under the equivalence of {\sl \'Etale\/} with presheaves on $L(S)$,
the representable presheaves correspond to
the \'etale actions $eS\ra E$, $s\mto s^*s= \mathbf{d}(s)$,
and the Yoneda embedding $L(S)\lra\PSh(L(S))$ is identified with the functor
\[
L(S)\lra \mbox{\sl \'Etale}\;;\; e\mto eS\ra E\;.
\]
A morphism $d\arr{s}e$ goes to the map
$\alpha_s\colon  dS\ra eS$ (over $E$) such that $\alpha_s(t)=st$.
For instance, $\alpha_s(d)=s$.
The Yoneda Lemma asserts in this case that $s\mto \alpha_s$
is a natural bijection between the \'{e}tale morphisms
$dS \ra eS$ and $L(S)(d,e)$.
Alternatively, we know that $C(S)(d,e)= eSd$ can be identified with morphisms $dS\ra eS$.
It is straightforward to verify that $s\in eSd$ corresponds
to a morphism over $E$ if and only if $s^*s=d$, i.e., $(e,s)\in L(S)$.

We proved in Lemma \ref{inpro} that the $S$-sets $eS=U(eS\ra E)$
are precisely the indecomposable projectives in
$S\text{-}\Se\/$ up to isomorphism.
Moreover, the functor $e\mto eS$ of $C(S)$ into $S\text{-}\Se\/$
is full and faithful, so that $C(S)$ is therefore weakly equivalent
to the full subcategory of $S\text{-}\Se\/$ on the indecomposable projectives.
When this functor is restricted to the subcategory $L(S)$,
the following diagram of functors commutes.
\begin{equation}
\xy
\square(0,0)|alra|/>`>`>`>/<600,400>%
  [L(S)`C(S)`\mbox{\sl \'Etale}`S\text{-}\Se;`{\rm Yoneda}``U]
\endxy
\label{Ufunctor}
\end{equation}
The functor $U(X\ra E)=X$
that forgets \'etale structure is faithful.

\begin{lemma}\named{yoneda}
Let $S$ be an inverse semigroup.
\begin{enumerate}
\item
A morphism of {\sl \'Etale\/}  is an epimorphism
if and only if it is a surjection.
\item
A morphism of {\sl \'Etale\/}
is a monomorphism if and only if it is injective.
In particular,  an \'{e}tale morphism $dS \ra eS$ is injective.
\end{enumerate}
\end{lemma}
\proof
The presheaf on $L(S)$ that corresponds to $X\arr{p}E$ is
the `fiber map' $e\mto p^{-1}(e)$.
If $d\arr{s}e$ in $L(S)$,
then the transition map for the presheaf
moves $x\in p^{-1}(e)$ to $xs\in  p^{-1}(d)$.
A morphism of \'etale actions is an epimorphism
if and only if its corresponding natural transformation of presheaves
is an epimorphism if and only if its component maps are surjections
if and only if the given map of \'etale actions is a surjection.
Alternatively, one can verify directly that a morphism
of {\sl \'Etale\/} is a epimorphism
if and only if the corresponding morphism of
$S\text{-}\Se\/$ is one,
and then use the corresponding result for $S\text{-}\Se\/$.
Both arguments can be repeated for monomorphisms and injections.

From a semigroup point of view,
a map $dS\arr{\alpha}eS$ (over $E$)
between representables is injective because
such a map is given by left multiplication
by an element $s\in eSd$ with $s^*s=d$: $\alpha(t)=st$.  The fact that multiplication on the left by $s$ is injective on $s^*sS$ is the trivial part of the classical Preston-Wagner theorem.
\qed

The \'etale version of Lemma \ref{inpro} is the following.

\begin{lemma}\named{projind}
An \'etale action $X\ra E$ is isomorphic to a representable one
$dS\ra E$ if and only if it is projective and indecomposable.
The Yoneda functor {\em (explained above)\/} gives a weak equivalence
between $L(S)$ and the full subcategory of {\sl \'Etale\/}
on the projective indecomposable objects.
\end{lemma}
\proof
This is a consequence of Prop.\ \ref{projconn}.
\qed

In the proof of Prop.\ \ref{CGCH} we encountered the fact that $C(S)$
is equivalent to ${\rm Span}(L(S))$.
Indeed, two functors
\[
\xy
\morphism(0,0)|a|/{@{>}@/^1em/}/<600,0>[C(S)`{\rm Span}(L(S));]
\morphism(0,0)|b|/{@{<-}@/_1em/}/<600,0>[\phantom{XX}`\phantom{XXXXX};]
\endxy
\]
giving the equivalence are $(e,s,d)\mto ((e,s),(d,s^*s))$,
and $((e,s),(d,t))\mto (e,st^*,d)$.
In terms of $S$-sets and \'etale actions,
an $S$-equivariant map $dS\arr{\theta} eS$ of $S$-sets
corresponds to a span of \'etale maps
\[
\xy
\Atriangle(0,0)|lrb|/>`>`/<300,300>[s^*sS`dS`eS;\theta_1`\theta_2`]
\endxy
\]
defined as follows:
$\theta_{1} (t) = ss^*t$, and $\theta_{2}(t) = st$.
Observe that $\theta_1$ is subset inclusion since $s^*s \leq d$.
Spans are composed in an obvious manner by pullback.

We return to the faithful functor $U$
that forgets \'etale structure (\ref{Ufunctor}).

\begin{proposition}
$U$ has a right adjoint $R$:
$$
R(X)=\coprod_E Xe \ra E\;;\;(e,x)\mto e\;,
$$
where
$$
Xe=\{x\in X\mid xe=x\}=\{xe \mid x\in X\}\iso S\text{-}\Se(eS,X)
$$
for an idempotent $e$.
For any $S$-set $X$, the counit $UR(X)\ra X$ is a surjection,
so that $R$ is faithful.
\end{proposition}
\proof
We denote a typical member of the coproduct
$\coprod_E Xe$ by $(e,x)$.
$\coprod_E eX$ is the sub-$S$-set of
$E\times X$ consisting of all pairs $\{(e,x)\mid xe=x\}$.
The action by $S$ that $\coprod_E eX$ carries is defined by:
$$
(e,x)s= (s^*es,xs)\;.
$$
Since idempotents commute in $S$,
if $e$ fixes $x$, then $s^*es$ fixes $xs$:  $xs(s^*es)=xess^*s=xs$.
The projection to $E$ is easily seen to be \'etale.
The unit of $U\adj R$ at $X\arr{p}E$ is the following map
of \'etale $S$-sets.
\begin{equation}
\xy
\Vtriangle(0,0)|alr|/>`>`>/<300,300>[X`\coprod_EXe`E;x\mto (p(x),x)`p`]
\endxy
\label{unitUR}
\end{equation}
The counit $UR(X)\ra X$
is the map $\coprod_EXe\ra X$,
$(e,x)\mto x$.
We have seen in the proof of Proposition~\ref{inpro} that
unitary is equivalent to the condition
$$
\forall x\in X,\exists e\in E, xe=x\;,
$$
which holds if and only if $\coprod_E Xe\ra X$ is onto.
\qed
$R$ may also be described as the equalizer:
\[
\xy
\morphism(0,0)|a|/ >->/<500,0>[R(X)`\phantom{XXX};]
\morphism(500,0)|a|/{@{>}@/^1em/}/<500,0>[E\times X`X;xe]
\morphism(500,0)|b|/{@{>}@/_1em/}/<500,0>[\phantom{XXXX}`\phantom{X};x]
\endxy\;.
\]
Evidently, $R$ is the composite
\begin{equation}
\xy
\btriangle(0,0)|lra|/>`>`>/<600,400>%
  [S\text{-}\Se`S\text{-}\Se/E`\mbox{\sl \'Etale};E^*`R`V]
\endxy
\label{VER}
\end{equation}
of two right adjoints,
where $E^*(X)=E\times X\ra E$,
and
\[
V(X\arr{p}E)=\{x\mid xp(x)=x\}\ra E\;,
\]
which is right adjoint to inclusion.
Because idempotents commute in $S$,
the action of $S$ in $X$ restricts to
$\{x\mid xp(x)=x\}$:
$$
xsp(xs)=xss^*p(x)s=xp(x)ss^*s=xs\;.
$$

\begin{lemma}
$R$ reflects isomorphisms.
\end{lemma}
\proof
Suppose that $X\arr{\psi}Y$ is a map of $S$-sets,
and that $R(\psi)$ is an isomorphism.
Then $\psi$ is a surjection because the counits of $U\adj R$ are surjections.
Now we prove that $\psi$ is injective.
Since $R(\psi)$ is injective, the restriction of $\psi$
to every $Xe$ is injective.
Suppose that $\psi(x)=\psi(x')$.
There are idempotents $d,e$ such that $xd=x$ and $x'e=x'$.
Then $\psi(xe)=\psi(x)e=\psi(x')$.
Since $xe,x'\in Xe$, we have $xe=x'$ by hypothesis.
Then $x'd=xed=xde=xe=x'$, so that $x,x'\in dX$.
Hence, $x=x'$ again since the restriction of $\psi$ to $Xd$ is injective.
Thus, $\psi$ is a bijection so that
it is an isomorphism in $S\text{-}\Se\/$.
\qed

Recall that  a {\em monad} \cite{BW} in a category
is an endofunctor $M$ of the category
equipped with natural transformations $M^2\ra M$ and ${\id}\ra M$,
called the multiplication and unit, respectively.
Associativity and unit conditions are required.
The (Eilenberg-Moore)
algebras for a monad form a category that maps to the given category
by forgetting an algebra's $M$ structure.
A functor is said to be {\em monadic} if it is equivalent to such a forgetful functor from the category of algebras for a monad.
We will use the following weak version of Beck's theorem:
if a functor has a left adjoint, reflects isomorphisms, coequalizers exist
and the functor preserves them, then it is monadic.
A {\em comonad} is a monad in the opposite category.
For all topos terminology and facts that we use, see \cite{MM}.

We begin by examining the restriction of presheaves
along the inclusion functor $I\colon L(S)\ra C(S)$, which we denote
$$
I^*\colon \PSh(C(S))\lra\PSh(L(S))\;.
$$
Under the equivalence of $\PSh(L(S))$ and {\sl \'Etale\/},
if $P$ is a presheaf on $C(S)$, then $I^*(P)$ is the \'etale action
$$
\coprod_E P(e)\ra E\;,
$$
where $(e,x)s=(s^*es,P(es)(x))$.
$I^*$ is the inverse image functor of a geometric morphism of toposes
$$
I^*\adj I_*\colon \mbox{\sl \'Etale}\lra\PSh(C(S))\;.
$$
The right adjoint $I_*$ is given by `taking sections,'
whose explicit description we omit.
The above geometric morphism is commonly termed a surjection because
its inverse image functor $I^*$ reflects isomorphisms.
Thus, in a geometric sense, $C(S)$ is a quotient of $L(S)$.
By the (dual) weak form of Beck's theorem,
$I^*$ is comonadic by a finite limit preserving comonad.
(A well-known fact of topos theory
is that a functor is equivalent to the inverse image functor
of a surjective geometric morphism if and only if it is comonadic
by a finite limit preserving comonad.)

$I^*$ also has a left adjoint $I_!$.
By definition,
if $X\arr{p}E$ is \'etale, and $e$ is an idempotent, then
\begin{equation}
I_!(p)(e) =
\colimit{\bX}\left(\; x\mto C(S)(e,p(x))\; \right)\;,
\label{Ipe}
\end{equation}
where $\bX$ is the category whose objects are the elements of $X$,
and morphisms $x\arr{s}y$ are
morphisms $p(x)\arr{s}p(y)$ of $L(S)$
satisfying $ys=x$.
$I^*$ is also monadic: it reflects isomorphisms, has a left adjoint,
and preserves all coequalizers.
The monad  $I^*I_!$ in {\sl \'Etale\/} associated with $I^*$
preserves all colimits,
and its category of algebras is equivalent to $\PSh(C(S))$.

Consider the following commuting diagrams of functors.
\[
\xy
\Atriangle(0,400)|lrb|/>`>`/<400,400>%
  [C(S)`\phantom{XXX}`\phantom{XXX}; F(e)= eS`{\rm Yoneda}`]
\square(0,0)|alra|/>`>`>`>/<800,400>%
  [S\text{-}\Se`\PSh(C(S))`S\text{-}\Se/E`%
      \mbox{\sl \'Etale};Q`E^*`I^*`V]
\btriangle(0,0)|lrb|/`>`/<700,400>%
  [\phantom{XXX}`\phantom{XXX}`\phantom{XXX};`R`]
\Atriangle(1400,400)|lrb|/>`>`/<400,400>%
  [C(S)`\phantom{XXX}`\phantom{XXX};F(e)=eS`{\rm Yoneda}`]
\square(1400,0)|alra|/<-`<-`<-`<-/<800,400>%
  [S\text{-}\Se`\PSh(C(S))`S\text{-}\Se/E`\mbox{\sl \'Etale};Q_!``I_!`]
\btriangle(1400,0)|lrb|/`<-`/<700,400>%
  [\phantom{XXX}`\phantom{XXX}`\phantom{XXX};`U`]
\endxy
\]
We have already met the functor $Q$ given by $Q(X)(e)= Xe$
and its left adjoint $Q_!$ in the proof of Theorem \ref{directproof}:
\[
Q_!(P)=\colimit{}\left( \bP\ra C(S)
  \arrow{F}S\text{-}\Se\right)\;,
\]
where $\bP\ra C(S)$ is the discrete fibration of elements of $P$.
$Q$ is faithful since $R$ is.
$I^*$ and $E^*$ are also faithful.
Of course,
the corresponding diagram of left adjoints commutes (above, right):
we have $Q_! I_!\iso U$, and $Q_!$ commutes with Yoneda.

\begin{lemma}\named{ISigU}
We have $I_!\iso Q U$:
for any \'etale $X\arr{p}E$ and any $e\in E$, $I_!(p)(e)\iso Xe$.
\end{lemma}
\proof
We argue this fact by direct calculation.
Let $X\arr{p}E$ be an \'etale action.
We claim that the unit map $I_!(p)\ra QQ_! I_!(p)\iso QU(p)$ is
a natural isomorphism (of presheaves on $C(S)$).
For any $e\in E$, the component map  at $e$ of this unit is
\[
I_!(p)(e)=\coprod_{x\in X} C(S)(e,p(x))/{\sim} \ra Xe\;;\;
\mbox{\rm equiv.\ class of}\;(x,e\arr{s}p(x))\mto xs\;,
\]
where the left-hand side is the colimit (\ref{Ipe}),
calculated as a coproduct factored by an equivalence relation.
This map has inverse $x\mto (x,e\arr{p(x)}p(x))$,
where $e\arr{p(x)}p(x)$ is the inequality
$p(x)\leq e$ understood as a map in $C(S)$,
which holds because $xe=x$, hence $p(x)e=p(x)$.
Furthermore, given any $(x,e\arr{s}p(x))$,
the map $xs\arr{s}x$ in the category $\bX$ (from \ref{Ipe})
witnesses that $(x,e\arr{s}p(x))$ is equivalent in the colimit
to $(xs, e\arr{p(xs)} p(xs))$,
noting $p(xs)=s^*p(x)s=s^*s\leq e$.
\qed

\begin{proposition}\named{UR}
$U$ reflects isomorphisms, $U$ has a right adjoint,
and {\sl \'Etale\/} has all equalizers
and $U$ preserves them.
$U$ is therefore comonadic.
\end{proposition}
\proof
$U$ preserves equalizers because they are created in both
categories by their underlying sets.
\qed

\begin{proposition}\named{Icomonad}
$I_!$ reflects isomorphisms, $I_!$ has a right adjoint,
and {\sl \'Etale\/} has all equalizers
and $I_!$ preserves them.
$I_!$ is therefore comonadic.
\end{proposition}
\proof
$I_!$ reflects isomorphisms because $U$ does
and $Q_! I_!\iso U$.
By Lemma \ref{ISigU},
$I_!$ preserves any limit $U$ does, such as an equalizer,
because $Q$ preserves all limits.
\qed

We have seen that $I^*$, $I_!$ and $U$ are all comonadic,
and that $I^*$ is also monadic,
but we wish to emphasize the following fact.

\begin{theorem}\named{EMUR}
$R$ is monadic.
The endofunctor of this monad
carries $X\arr{p}E$ to $\coprod_E Xe\ra E$,
as in (\ref{unitUR}).
In other words,
its category of Eilenberg-Moore algebras is equivalent to $S$-\Se.
\end{theorem}
\proof
To show that $R$ is monadic it suffices to show that $R$
preserves coequalizers since we already know that $R$
reflects isomorphisms and has a left adjoint.
We shall do this by inspecting the construction of coequalizers,
which is relatively straightforward since coequalizers are set-theoretic
in both {\sl \'Etale\/} and $S$-\Se.
Let
\[
\xy
\morphism(500,0)|a|/>>/<500,0>[\phantom{X}`C;\psi]
\morphism(0,0)|a|/{@{>}@/^1em/}/<500,0>[X`Y;f]
\morphism(0,0)|b|/{@{>}@/_1em/}/<500,0>[\phantom{X}`\phantom{X};g]
\endxy
\]
be a coequalizer in $S$-\Se.
Applying $R$ gives a diagram
\[
\xy
\morphism(500,0)|a|/>>/<500,0>[\phantom{XXX}`K;]
\morphism(1000,0)|b|/>/<500,0>[\phantom{X}`\coprod_E Ce;\eta]
\morphism(0,0)|a|/{@{>}@/^1em/}/<500,0>[\coprod_E Xe`\coprod_E Ye;]
\morphism(0,0)|b|/{@{>}@/_1em/}/<500,0>[\phantom{XXX}`\phantom{XXX};]
\morphism(500,0)|a|/{@{>}@/^1em/}/<900,0>[\phantom{XXX}`\phantom{XX};R(\psi)]
\endxy
\]
where $K$ is the coequalizer in {\sl \'Etale\/}.
$R(\psi)$ is a surjection since given $c\in Ce$,
there is $y\in Y$ such that $\psi(y)=c$.
Then $\psi(ye)=\psi(y)e=ce=c$, and $ye\in Ye$.
Therefore, $\eta$ is a surjection.
$\eta$ is also injective: suppose that
$R(\psi)(d,y)=R(\psi)(e,y')$.
Then $d=e$ and $\psi(y)=\psi(y')$.
This says that $y$ and $y'$ are connected by a
finite `zig-zag' under $f$ and $g$.
For instance, we may have a two-step zig-zag
\[
\xy
\Atriangle(0,0)|lrb|/>`>`/<300,300>[x`y`y'';f`g`]
\Atriangle(600,0)|lrb|/>`>`/<300,300>[x'`\phantom{XX}`y';f`g`]
\endxy
\]
Multiply the zig-zag by $d$ so
that $(d,y)$ and $(d,y')$ are equal in $K$.
This shows that $\eta$ is injective,
whence an isomorphism in {\sl \'Etale}.
\qed

We may now deduce the inverse case of Theorem \ref{directproof}
in a different way.

\begin{corollary}\named{SCS}
The monads in {\sl \'Etale\/}
associated with the adjoint pairs $U\adj R$
and $I_!\adj I^*$ are isomorphic.
{\em (Thus, this monad preserves all colimits.)}
The adjoint pair
\begin{equation}
Q_!\adj Q\colon S\text{-}\Se\equ\PSh(C(S))
\label{SsetCS}
\end{equation}
is an equivalence.
\end{corollary}
\proof
The two monads $RU$ and $I^*I_!$ are isomorphic because,
by Lemma \ref{ISigU}, we have $I^*I_!\iso I^*QU\iso RU$.
The two monads therefore have equivalent algebra categories:
for $I^*I_!$ it is $\PSh(C(S))$, and for $RU$ it is $S\text{-}\Se\/$
(Thm.\ \ref{EMUR}).
\qed

The fact that $\PSh(C(S))$
and $S\text{-}\Se\/$ are equivalent
generalizes the well-known fact when $S=M$ is an (inverse) monoid that
presheaves on a category and on its Cauchy completion are equivalent
because $C(M)$ is the Cauchy completion of
$M$ as a category (with a single object)~\cite{Bo}.


\section{Complements}\named{C}

There is a variation of enlargement Morita equivalence that uses only semigroups.
However, the Axiom of Choice is used.
Lawson \cite{L1} generalized the property of an idempotent $e$ that $S = SeS$.
If $S$ is a subsemigroup of another semigroup $T$ we say that
$T$ is an {\em enlargement\/} of $S$ if $S = STS$ and $T = TST$.
If $S = SeS$, then $S$ is an enlargement of $eSe$.
Lawson \cite{L3} observes that if $S$ and $T$ have local units
and $T$ is an enlargement of $S$,
then $S$ and $T$ are Morita equivalent in the Talwar sense.
If $R$ is an enlargement of subsemigroups $S$ and $T$,
then we say that $R$ is a {\em joint enlargement} of $S$ and $T$.
If $R$ is a regular, then we say that it is a {\em regular} joint enlargement.

\begin{theorem}[Axiom of Choice]\label{renlar}
Inverse semigroups $S$ and $T$ are  strongly Morita equivalent
if and only if there is a regular semigroup
that is a joint enlargement of $S$ and $T$.
\end{theorem}
\proof
If $S$ and $T$ are strongly Morita equivalent, then $C(S)$ and $C(T)$ form a Morita context by Proposition~\ref{CGCH}.
Lawson \cite{AL} has proved in a more general frame that
this implies that $S$ and $T$ have a regular joint enlargement.

Conversely,
let the regular semigroup $R$ be a joint enlargement of
inverse subsemigroups $S$ and $T$.
Let $x \in SRT$.
Then $x = srt$.
Let $s^{\ast}$ be the unique inverse of $s$ in $S$,
and let $t^{\ast}$ be the unique inverse of $t$ in $T$.
Then $x$ has an inverse of the form $t^{\ast}r's^{\ast} \in TRS$,
where $r' \in R$ is some element.
Put
$$
X = \{(x,x') \colon \: x \in SRT \text{ and } x' \in V(x) \cap TRS \}.
$$
Observe that
$$
xx' \in (SRT)(TRS) = S(RTTR)S \subseteq S
$$
and
$$
x'x \in (TRS)(SRT) = T(RSSR)T \subseteq T\;.
$$
Thus we may define a left action of $S$ on $X$
by $s(x,x') = (sx,x's^{\ast})$
and a right action of $T$ on $X$ by $(x,x')t = (xt,t^{\ast}x')$.
Thus $X$ is an $(S,T)$-biset.
Define $\langle (x,x'),(y,y') \rangle = xy'$ and $[(x,x'),(y,y')] = x'y$.
We need to show that these maps are surjections.
We prove that the first is surjective;
the proof that the second is surjective follows by symmetry.
Let $s \in S$.
Then $s = bta'$ where $aa' = s^{\ast}s$
and $bb' = ss^{\ast}$, and $a \in V(a)$ and $b \in V(b)$.
A proof that this is possible is given in \cite{L1}.
Let $t \in V(t)$ such that $t't = a'a$ and $tt' = b'b$.
Then $(b,b'),(at',ta') \in X$ and
$\langle (b,b'), (at',ta') \rangle = bta' = s$, as required.
It is now routine to verify that axioms (M1) - (M7) hold
and that we have therefore defined an equivalence biset.
\qed

\begin{remark}{\em
The above result raises the following question:
is it true that two inverse semigroups which are Morita equivalent have a joint {\em inverse} enlargement?
We suspect this is not true, although we do not have a counterexample.
However, in the light of Proposition~5.9 \cite{S}
we make the following conjecture.
We say that an inverse semigroup $S$
is {\em directed} if for each pair of idempotents $e,f \in S$ there is an idempotent $i$ such that $e,f \leq i$.
This is equivalent to the condition that each subset of the
form $eSf$ is a subset of some local submonoid $iSi$.
Semigroups with this property are studied in \cite{N1,N2}.
We conjecture that if $S$ and $T$ are both directed, then they are Morita equivalent if and only if they have an inverse semigroup joint enlargement.}
\end{remark}


\begin{remark}
{\em If two inverse semigroups $S$ and $T$ have a
regular semigroup as a joint enlargement,
then it is easy to show that $C(S)$ and $C(T)$ are part of
a Morita context so that $S$ and $T$ are strongly Morita equivalent.
This does not require the Axiom of Choice.
However, we currently know of no proof of the converse
that does not use the Axiom of Choice.}
\end{remark}

We include here a direct proof that strong Morita equivalence
and enlargement equivalence are the same.
It uses the fact that we may generalize semigroups to {\em semigroupoids},
which are categories possibly without identities, but with objects.
Thus, a semigroup is a semigroupoid with one object.

\begin{theorem}\named{boge}
Two inverse semigroups are strongly Morita equivalent if and only if their associated inductive groupoids have a bipartite ordered groupoid enlargement.
\end{theorem}
\proof
Let $(S,T,X,\langle -,- \rangle,[-,-])$ be an equivalence biset.
Put $I = \{1,2 \}$ and regard $I \times I$ as a groupoid in the usual way,
$S' = \{1 \} \times S \times \{1\}$ and $T' = \{2 \} \times T \times \{2\}$
and
$$
\mathcal{R} = \mathcal{R}(S,T;X)
= S'
\cup
T'
\cup
(\{1 \} \times X \times \{2\})
\cup
(\{2 \} \times X \times \{1\})\;.
$$
We shall define a partial binary operation on $\mathcal{R}$.
The product of $(i,\alpha,j)$ and $(k,\beta, l)$ will be defined
if and only if $j = k$
in which case the product will be of the form $(i,\gamma,l)$.
Specifically, we define products as follows
\begin{itemize}

\item $(1,s,1)(1,s',1) = (1,ss',1)$.

\item $(2,t,2)(2,t',2) = (2,tt',2)$.

\item $(1,s,1)(1,x,2) = (1,sx,2)$.

\item $(1,x,2)(2,t,2) = (1,xt,2)$.

\item $(2,t,2)(2,x,1) = (2,xt^{\ast},1)$.

\item $(2,x,1)(1,s,1) = (2,s^{\ast}x,1)$.

\item $(2,x,1)(1,y,2) = (2,[x,y],2)$.

\item $(1,x,2)(2,y,1) = (1,\langle x,y\rangle,1)$.
\end{itemize}
This operation is associative whenever it is defined.
To prove this one essentially checks all possible
cases of triples of elements;
however, the restrictions on what elements
can be multiplied reduces the number of cases that need to be checked.
Within this list of possibilities,
associativity of multiplication in the inverse semigroups $S$ and $T$
combined with the `associativity' of left,
right and biset actions reduces the number of cases still further.
One then uses the definition of an equivalence biset, and particularly Proposition~2.3 of \cite{S}, to check all the remaining cases.
Thus $\mathcal{R}$ is a semigroupoid.
Observe that $(1,x,2)(2,x,1) = (1,\langle x,x \rangle, 1)$
and that $(2,x,1)(1,x,2) = (2,[x,x],2)$.
Thus
$$
(1,x,2)(2,x,1)(1,x,2) = (1,\langle x, x \rangle x  ,2) = (1,x,2)
$$
by (M3).
Similarly
$$
(2,x,1)(1,x,2)(2,x,1) = (2, [x,x] ,2)(2,x,1) = (2,x[x,x],1) = (2,x,1)
$$
by (M6).
Thus $\mathcal{R}$ is a regular semigroupoid.
But the only idempotents in $\mathcal{R}$ are
those coming from $S'$ and $T'$,
so that idempotents commute whenever
the product of two idempotents is defined.
It follows that $\mathcal{R}$ is an inverse semigroupoid.
Clearly $S' = S'\mathcal{R}S'$ and $T' = T'\mathcal{R}T'$,
and it is easy to check that
$\mathcal{R} = \mathcal{R}S'\mathcal{R}$ and $\mathcal{R} = \mathcal{R}T'\mathcal{R}$.
Every inverse semigroupoid gives rise to an ordered groupoid in a way that directly generalizes the way in which inverse semigroups
give rise to ordered groupoids.
We denote this ordered groupoid by
\begin{equation}
G(S,T;X)\;.
\label{GSTX}
\end{equation}
We see that $G(S,T;X)$ is an enlargement of both $G(S')$ and $G(T')$.

Conversely, let $S$ and $T$ be inductive groupoids which are
ordered subgroupoids of the ordered groupoid $G$,
and where $G$ is an enlargement of them both.
Let $X$ be the set of all the arrow of $G$
that have domains in $T$ and codomains in $S$.
We define a left action of $S$ on $X$ by $sx = s \circ x$,
and a right action of $T$ on $X$ by $xt = x \circ t$.
Define $\langle x,y \rangle = x \circ y^{-1}$,
and $[x,y] = x^{-1} \circ y$.
Here $\circ$ is the pseudoproduct in the ordered groupoid $G$.
It is routine using the theory of ordered groupoids and pseudogroups \cite{L2}
to check that in this way we have defined an equivalence biset.
\qed

We conclude this section with an application of Morita equivalence
to the theory of $E$-unitary inverse semigroups.
With each $E$-unitary inverse semigroup $S$
we can associate a triple $(G,X,Y)$, called a {\em McAlister triple},
where $G$ is a group, $X$ a poset, and $Y$ a downset of $X$
that is a semilattice for the induced order \cite{L2}.
This triple is required to satisfy certain conditions,
one of which is that $G$ acts on $X$ by order automorphisms.
If $(G,X)$ and $(G',X')$ each consist of a group acting
by order automorphisms on a poset,
then we say they are {\em equivalent\/} if there is a group isomorphism $\phi\colon G\ra G'$ and an order-isomorphism $\theta\colon X\ra X'$ such that $\theta (xg)=\theta(x)\phi(g)$ for all $x\in X$ and $g\in G$.

\begin{proposition}
Let $S$ and $T$ be $E$-unitary inverse semigroups with associated McAlister triples $(G,X,Y)$ and $(G',X',Y')$.
Then $S$ and $T$ are Morita equivalent if and only if
$(G,X)$ is equivalent to $(G',X')$.
\end{proposition}
\proof
Let $S$ and $T$ be such that $(G,X)$ is equivalent to $(G',X')$.
Then after making appropriate identifications,
we have from the classical theory of $E$-unitary
inverse semigroups \cite{L2} that
the Grothendieck or semidirect product construction $G\ltimes X$,
which is an ordered groupoid,
is a common enlargement of the inductive groupoids $G(S)$ and $G(T)$.

Conversely, suppose that $S$ and $T$ are strongly Morita equivalent.
Then the toposes $\B(S)$ and $\B(T)$ are equivalent.
The topos explanation of the $P$-theorem is simply
an interpretation of $X$ and $G$ in topos terms~\cite{F1}:
the (connected) universal covering morphism
of the classifying topos $\B(S)$ has the form $\PSh(X)\ra \B(S)$,
$G$ is the fundamental group of $\B(S)$,
and the action of $G$ on $X$ is induced from the action by deck transformations.  So if $\B(S)$ and $\B(T)$ are equivalent toposes,
then $(G,X)$ and $(G',X')$ must be equivalent.
An explicit description of an equivalence of $(G,X)$ and $(G',X')$
derived directly from and in terms of a given equivalence biset
ought to be readily available, but we leave this exercise for the reader.
\qed

Let us say that an  inverse semigroup $S$ is
{\em locally $E$-unitary\/}
if the local submonoid $eSe$ is $E$-unitary for every idempotent $e$.
An $E$-unitary inverse semigroup is locally $E$-unitary.

\begin{lemma}\named{rcanc}
$S$ is locally $E$-unitary if and only if $L(S)$ is right-cancellative.
\end{lemma}
\proof
Suppose that $L(S)$ is right-cancellative.
Let $s=ese$ and suppose that $d\leq s$, where $d$ is an idempotent.
Then the diagram $d\leq s^*s\arr{s,s^*s}e$ in $L(S)$ commutes.
Therefore, $s=s^*s$ so that $s$ is an idempotent.

Conversely, suppose that $S$ is locally $E$-unitary.
Suppose that $d\arr{t}e\arr{s,r}f$ commutes in $L(S)$.
Then $rs^*\in fSf$.
Also $rtt^*s^*=rt(st)^*=st(st)^*$ is idempotent,
and we have $rtt^*s^*\leq rs^*$.
Therefore, $rs^*=b$ is an idempotent by locally $E$-unitary.
Hence, $r=rr^*r=re=rs^*s=bs$, so that $r\leq s$.
Similarly, $s\leq r$ so that $s=r$.
\qed
We take the opportunity to improve \cite{F1}, Cor.\ 4.3.

\begin{corollary}
$\B(S)$ is  locally decidable {\em (as it is called)\/} if and only if
$S$ is locally $E$-unitary.
\end{corollary}
\proof
This follows from Lemma \ref{rcanc} and the well-known fact
that the topos of presheaves on a small category
is locally decidable if and only if
the category is right-cancellative \cite{El}.
\qed

\begin{corollary}
If two inverse semigroups are Morita equivalent
and one of them is locally $E$-unitary, then so is the other one.
\end{corollary}

{\small \bibliographystyle{plain}
}


\bigskip
jonathon.funk@cavehill.uwi.edu

The University of the West Indies, Bridgetown, Barbados.

\medskip
M.V.Lawson@ma.hw.ac.uk

Maxwell Institute for Mathematical Sciences, Mathematics Department,

Heriot-Watt University, Edinburgh, Scotland.

\medskip
bsteinbg@math.carleton.ca

Carleton University, Ottawa, Canada.

\end{document}